\documentclass[11pt]{article}
\usepackage{amscd, amsmath, amssymb}

\title{The resolution property for schemes and stacks}
\author{Burt Totaro}
\date{ }

\def\Z{\text{\bf Z}}
\def\Q{\text{\bf Q}}
\def\R{\text{\bf R}}
\def\C{\text{\bf C}}
\def\P{\text{\bf P}}
\def\E{{\cal E}}
\def\G{{\cal G}}
\def\Spec{\text{Spec }}
\def\Rep{\text{Rep}}
\def\Aff{\text{Aff}}

\def\arrow{\rightarrow}

\def\imp{\Rightarrow}
\def\surj{\twoheadrightarrow}
\def\qed{\ QED \vspace{0.15in}}
\def\naive{\text{naive}}

\begin{document}
\maketitle

\newtheorem{theorem}{Theorem}[section]
\newtheorem{corollary}[theorem]{Corollary}
\newtheorem{lemma}[theorem]{Lemma}
\newtheorem{proposition}[theorem]{Proposition}

{\bf Abstract. }We prove the equivalence of two fundamental properties
of algebraic stacks: being a quotient stack in a strong sense, and
the resolution property, which says that every coherent sheaf is
a quotient of some vector bundle. Moreover, we prove these properties
in the important special case of orbifolds whose associated
algebraic space is a scheme. (Mathematics Subject Classification:
Primary 14A20, Secondary 14L30.)

\section{Introduction}

Roughly speaking, an algebraic stack is an object which looks locally like
the quotient of an algebraic variety by a group action \cite{LMB}. Thus it
is
a fundamental question whether a given stack is globally the quotient
of a variety by a group action. We show that a strong version of
this property is equivalent to another fundamental property
of stacks: having ``enough'' vector bundles for geometric purposes.
The equivalence turns out to be interesting even in
the special case where the stack is a scheme. Moreover, we prove
the existence of enough vector bundles in the important case
of orbifolds whose associated algebraic space is a scheme.
In more detail, the two main results of the paper are:

\begin{theorem}
\label{main}
Let $X$ be a normal noetherian algebraic stack (over $\Z$) whose
stabilizer groups at closed points of $X$ are affine.
The following are equivalent.

(1) $X$ has the resolution property: every coherent sheaf on $X$
is a quotient of a vector bundle on $X$.

(2) $X$ is isomorphic to the quotient stack
of some quasi-affine scheme by an action of the group $GL(n)$
for some $n$.

For $X$ of finite type over a field $k$, these are also equivalent to:

(3) $X$ is isomorphic to the quotient stack of some affine scheme over $k$
by
an action of an affine group scheme of finite type over $k$.
\end{theorem}

\begin{theorem}
\label{orb2}
Let $X$ be a smooth Deligne-Mumford stack over a field $k$. Suppose that
$X$ has finite stabilizer group and that the stabilizer group is generically
trivial. Let $B$ be the Keel-Mori coarse moduli space of $X$ \cite{KM}. If
the algebraic space
$B$ is a scheme with affine diagonal (for example a separated scheme),
then the stack 
$X$ has the resolution property.
\end{theorem}

Informally, Theorem \ref{orb2}
 says that any ``orbifold coherent sheaf'' on a scheme with
quotient singularities admits a resolution by ``orbifold vector bundles.''
For example, Theorem \ref{orb2} implies
that the moduli stacks of curves, $\overline{M}_{g,n}$ with
$g\geq 3$ (so that the generic stabilizer is trivial), have the
resolution property. This was proved earlier by Mumford \cite{Mumfordcurve}
by showing that these particular stacks admit a Cohen-Macaulay global
cover. The more general Theorem \ref{orb2} has often been wished for,
even in the special case where $B$ is quasi-projective.
It allows the
hard-to-verify assumption of a Cohen-Macaulay global cover to be removed
from various papers, such as Kawamata's paper on flips and derived
categories \cite{Kawamata}.

In Theorem \ref{main},
the fact that quotient stacks $W/GL(n)$ with $W$ a noetherian
quasi-affine scheme (an open subset
of an affine scheme) have the resolution property is an easy special case of
Thomason's results on the resolution property \cite{Thomasonres},
listed in Theorem \ref{thom} below.
So the new
implication is the opposite one, which says that stacks with the resolution
property
are quotient stacks of a very special kind.
Edidin, Hassett, Kresch, and Vistoli proved a 
step in the direction of Theorem \ref{main}: they showed that a stack
with quasi-finite
stabilizer group which has the resolution property
is a quotient stack $W/GL(n)$ with $W$ an algebraic space, that is,
a stack with trivial stabilizer group
(\cite{EHKV}, Theorem 2.7 and Theorem 2.14). Algebraic spaces (and even
schemes) do not all have the resolution property, as explained below,
and so we need to
strengthen the conclusion as in Theorem \ref{main}
in order to have an equivalence.

{\bf Remarks. }(1)
The restriction to stacks with affine stabilizer groups in Theorem
\ref{main} seems reasonable.
In fact,
I would argue that the resolution property is not meaningful for a stack
whose
stabilizer
groups are not affine. The simplest example of such a stack is the
classifying stack $BE$
of an elliptic curve $E$. Because every linear representation of an elliptic
curve is trivial,
coherent sheaves and vector bundles on $BE$ are both simply vector spaces.
Thus, the resolution property for $BE$ holds, but this has no real geometric
significance. In particular, the $K$-groups of $BE$ defined
using either vector bundles or perfect complexes (cf.\ section
\ref{history})
are simply the
$K$-groups of a point.

(2) Theorem \ref{main} cannot be strengthened to say that a stack $X$ with
the resolution
property (and affine stabilizer groups) is a quotient stack
$W/GL(n)$ with $W$ an
affine scheme. Indeed, quotient stacks of the latter form are very special.
For example,
by geometric invariant theory, any algebraic space which is the quotient
of an affine scheme by a reductive group such as $GL(n)$
(as a stack, which means that the
group action is free) is in fact an affine scheme (\cite{GIT}, Amplification
1.3).

Theorem \ref{main} implies:

\begin{proposition}
\label{affine}
Let $X$ be a noetherian algebraic stack whose stabilizer groups at
closed points of $X$ are affine. If $X$ has
the resolution property,
then the diagonal morphism $X\arrow X\times_{\Z} X$ is affine.
\end{proposition}

For example, a scheme $X$ has affine diagonal morphism if and only if the
intersection
of any two affine open subsets is affine. The property of having affine
diagonal is a natural weakening of separatedness: the smooth non-separated
scheme
$A^n\cup_{A^n-0}A^n$ has affine diagonal for $n=1$ but not for $n\geq 2$.
Thus Proposition \ref{affine} immediately implies Thomason's observation that
the scheme
$A^n\cup_{A^n-0}A^n$ does not have the resolution property for $n\geq 2$
(\cite{TT}, Exercise 8.6). Most other
known counterexamples to the resolution property have non-affine
diagonal and so are explained by
Proposition \ref{affine}. Unfortunately, not all stacks with affine
diagonal have the resolution property. By Grothendieck, there is
a normal
(affine or projective) complex surface $Y$ which has a non-torsion
element of $H^2_{\text{et}}(Y,G_m)$
(\cite{GrothendieckBrauer}, II.1.11(b)). As Edidin, Hassett, Kresch,
and Vistoli observed, the corresponding
$G_m$-gerbe $X$ over $Y$ (a stack with stabilizer group $G_m$ at every point)
is a stack with affine diagonal which is not
a quotient of an algebraic space by $GL(n)$ (\cite{EHKV}, Example 3.12).
Therefore, Theorem \ref{main} says that $X$ fails to have the
resolution property. Despite this counterexample, we can
ask the following optimistic questions about the resolution property:

{\bf Question 1. }Let $X$ be a noetherian stack with quasi-finite
stabilizer group (for example, an algebraic space or a scheme).
Suppose that $X$ has affine diagonal (for example, $X$ separated).
Does $X$ have the resolution property?

This would be very useful. I do not know how likely it is.
At the moment, the resolution property is not known even for such concrete
objects as normal toric varieties or smooth separated algebraic spaces
over a field.

{\bf Question 2. }Let $X$ be a smooth stack over a field such that $X$
has affine
diagonal. Does $X$ have the resolution property?

{\bf Question 3. }Let $X$ be a stack which has the resolution property.
Suppose that $X$ has an action of a flat affine group scheme $G$
of finite type over
a field or over $\Z$. Does $X/G$ have the resolution property?

The converse to Question 3 is true and sometimes useful: if $G$
acts on a stack $X$ (for example, a scheme) such that the quotient
stack $X/G$ has the resolution
property, then so does $X$ (Corollary \ref{quotres}).

Thomason proved
several cases of Question 3, listed in Theorem \ref{thom} below.
One case where his methods do not immediately apply
is the quotient stack $Q$ of the nodal cubic curve by the multiplicative
group, which has several pathological properties described in
section \ref{nodal}. Nonetheless, in that section we use some ideas
from the proof of Theorem \ref{main} to prove the resolution
property for $Q$. In general, the proof
of Theorem \ref{main} often shows how to
construct a coherent sheaf $C$ on a given stack $X$ such that $X$ has
the resolution property if and only if the single sheaf $C$ is a quotient
of a vector bundle. I expect that this method should help to prove
the resolution property in other situations as well.

Finally, section \ref{counterexamples} considers the question of
whether surjectivity of the natural map from 
the Grothendieck group $K_0^{\naive}X$ of vector bundles on $X$ to the
group $K_0X$ of perfect complexes is enough to imply the resolution
property. The answer is yes for smooth schemes, but not for smooth
algebraic spaces or Deligne-Mumford stacks, as we will show in
two examples.

I would like to thank
Yujiro Kawamata and Gabriele Vezzosi
for useful discussions on Theorem \ref{orb2}. In particular, the proof
given here of Theorem \ref{orb2} has been simplified
by an idea of Vezzosi's.

\tableofcontents

\section{History of the resolution property}
\label{history}

One important reason to consider the resolution property is its role in
$K$-theory. For any scheme $X$, it is natural to consider the Grothendieck
group of vector bundles $K_0^{\naive}X$. Quillen extended this group to
a sequence of groups $K_*^{\naive}X$ built from the category of vector
bundles on $X$. We use the name ``naive'' because these groups
do not satisfy the important properties such as Mayer-Vietoris for
open coverings of arbitrary schemes. Thomason showed how to define
the right $K$-groups $K_*X$ for arbitrary schemes \cite{TT};
more generally, it seems that the same definition works for a
stack $X$.
Instead of vector bundles, he used perfect complexes, that is, complexes
of sheaves which are locally quasi-isomorphic to bounded complexes
of finite-dimensional vector bundles. The same idea comes up in topology.
In simple situations such as compact Hausdorff spaces,
one can define topological $K$-groups
using only finite-dimensional vector bundles. In more complex situations,
in order to get the right $K$-groups (satisfying Mayer-Vietoris),
one has to consider Fredholm complexes, that is, complexes of
infinite-dimensional vector bundles with finite-dimensional cohomology
sheaves. An example where infinite-dimensional bundles are needed
is the twisted $K$-theory recently considered by Freed, Hopkins,
Teleman, and others, which can be viewed as the $K$-theory of a
``topological stack,'' the $S^1$-gerbe $X$ over a space $Y$ associated
to an element of $H^2(Y,S^1_{\text{cont}})=H^3(Y,\Z)$; when this element is not
torsion, there are not enough finite-dimensional vector bundles on the
stack $X$ (\cite{Atiyahtwist}, \cite{FHT}).

It is natural to ask for criteria to ensure that
the natural map from $K_*^{\naive}X$ to $K_*X$
is an isomorphism,
because this means that Thomason's $K$-groups give information
about vector bundles on $X$, which are geometrically appealing.
Thomason gave a satisfactory answer: if $X$ has the resolution property,
then the map from
$K_*^{\naive}X$ to $K_*X$ is an isomorphism.

For clarity, perhaps I should
add that the $G$-theory (or $K'$-theory)
of coherent sheaves, also defined by Quillen, has good
properties in general and has not had to be modified. In particular, for
any regular scheme $X$, the natural map $K_*X\arrow G_*X$ is an isomorphism,
whereas these groups may differ from $K_*^{\naive}X$ when the resolution
property
fails, for example for the regular scheme
$X=A^n\cup_{A^n-0}A^n$ with $n\geq 2$ (\cite{TT}, Exercise 8.6).

The resolution property is known to hold for a vast class of schemes
and stacks. In particular, it holds for any noetherian scheme with an
ample family of line bundles, by Kleiman and independently
Illusie. Kleiman's proof is given in Borelli \cite{Borelli67},
3.3, and Illusie's is in \cite{Illusie}, 2.2.3 and 2.2.4. A convenient
summary of the results on ample families of line bundles is given
in Thomason-Trobaugh \cite{TT}, 2.1.
By definition, a scheme $X$
has an ample family of line bundles if $X$ is the union of open affine
subsets of the form $\{ f\neq 0\}$ with $f$ a section of a line bundle
on $X$.
Borelli and Illusie also showed that
every regular (or, more generally, factorial)
separated noetherian scheme has an ample family of line bundles
(\cite{Borelli63}, 4.2 and \cite{Illusie}, 2.2.7). (In fact,
the same proof works a little more generally, as observed recently by
Brenner and
Schr\"oer: every $\Q$-factorial noetherian scheme with affine diagonal has an
ample family of line bundles (\cite{BS}, 1.3).)
Also, an ample family of line bundles passes to locally closed subschemes,
and so the resolution property holds for all quasi-projective schemes
over an affine scheme.
Finally, Thomason proved the resolution property for the most
naturally occurring stacks, quotient stacks of the following types among
others
(\cite{Thomasonres}, 2.4, 2.6, 2.10, 2.14):

\begin{theorem}
\label{thom}
Let $X$ be a noetherian scheme over a regular noetherian ring $R$
of dimension at most 1.
Let $G$ be a flat affine group scheme of finite type over $R$
together with an action of $G$ on $X$.

(1) If $X$ has an ample family of $G$-equivariant line bundles,
then the stack $X/G$ has the
resolution property.

(2)  Suppose that $G$ is an extension of a finite flat group scheme
by a smooth group scheme with connected fibers over $R$. (This is
automatic if $R$ is a field.)
If $X$ is normal and has an ample family of line bundles,
then the stack $X/G$ has the resolution property.
\end{theorem}

Equivalently, under these
assumptions,
every $G$-equivariant coherent sheaf on $X$ is a quotient of a
$G$-equivariant
vector bundle on $X$. For example, it follows
that the resolution property holds for one of Hironaka's examples
of a smooth proper
algebraic space $X$ of dimension 3 which is not a scheme: $X$
is the quotient of a smooth separated scheme (hence a normal scheme with an
ample
family of line bundles) by a free action of
the group $\Z/2$ (\cite{Knutson}, pp.~15--17). Hironaka's paper
\cite{Hironaka}, Example 2, gives closely
related examples for which the resolution property seems to be unknown.

Hausen showed recently that a scheme (assumed to be reduced and of finite
type
over an algebraically closed field, although that should be unnecessary)
has an ample family of
line
bundles if and only if it is a quotient $W/(G_m)^n$ for some free action
of the torus $(G_m)^n$ on a quasi-affine scheme $W$ (\cite{Hausen}, 1.1).
This is a satisfying analogue to Theorem \ref{main}, though the situation
is simpler in that only schemes are involved.

We can ask whether the resolution property holds in cases not covered
by the above results, but here progress has been slow. As mentioned
in the introduction, the resolution property is still an open problem even
for such
concrete objects as normal toric varieties or smooth separated algebraic
spaces. For example, Fulton found a 3-dimensional normal
proper toric variety $X$ over $\C$ which has zero Picard group
(\cite{Fulton}, p.~65);
such a variety cannot have an ample family of line bundles,
and it is unknown whether the resolution property holds. One
encouraging result is the recent proof by Schr\"oer and Vezzosi of
the resolution property for all normal separated surfaces (\cite{SV}, 2.1).
There
is a normal proper surface over $\C$ that has zero Picard group.
Thus, Schr\"oer and Vezzosi were able to construct enough vector bundles
to prove the resolution property on such a variety, even though it
does not have an ample family of line bundles.

As indicated in the introduction, all known examples of
stacks without the resolution property are non-separated.
The situation is very different in the complex analytic
category, where Voisin has recently proved that the
resolution property can fail even for compact K\"ahler manifolds of
dimension 3 (\cite{Voisin}, Appendix).

\section{Proof of Theorem \ref{main}, first implications}
\label{easy}

We begin by proving the easier parts of Theorem \ref{main}: that (2) implies
(1) and that (2) and (3) are equivalent.

First, (2) implies
(1). Let $X$ be
the quotient stack $W/GL(n)$ for some action of the group $GL(n)$ (over
$\Z$)
on a
noetherian quasi-affine scheme $W$.
A noetherian scheme $W$ is quasi-affine if and only
if the trivial line bundle $O_W$ is ample (\cite{EGA}, Proposition
II.5.1.2). Given an action of
$GL(n)$
on $W$, the trivial line bundle $O_W$
has the structure of a $GL(n)$-equivariant line bundle in a natural way,
so $O_W$ is a $GL(n)$-equivariant ample line bundle on $W$. By 
Theorem \ref{thom}, due to Thomason,
it follows that every $GL(n)$-equivariant coherent sheaf
on $W$ is a quotient of some $GL(n)$-equivariant vector bundle.
Equivalently,
the stack $W/GL(n)$ has the resolution property.

Next, let us show that (2) implies (3) in Theorem \ref{main}. Here we only
consider stacks of finite type over a field $k$. The proof is an equivariant
version of Jouanolou's trick, which we will prove directly (\cite{Jouanolou},
1.5).
Let $X$ be a quotient stack $W/GL(n)$ for some action of $GL(n)$
over $k$ on a quasi-affine scheme $W$
of finite type over $k$. By EGA II.5.1.9 \cite{EGA}, $W$ embeds as an open
subset of an
affine scheme $Y$ of finite type over $k$; we can also arrange for this
embedding
to be $GL(n)$-equivariant, simply by taking $\Spec$of a bigger finitely
generated
subalgebra of $O(W)$. Here we are using that the $GL(n)$-module $O(W)$ over
$k$
is a union of finite-dimensional representations of $GL(n)$;
a direct proof
is given in GIT \cite{GIT}, p.~25,
although it is also a special case of the fact
that
a quasi-coherent sheaf on a noetherian stack (here $BGL(n)$ over $k$)
is a filtered direct limit of
its coherent subsheaves (\cite{LMB}, 15.4).
 Let the closed subset $Y-W$ be defined by the vanishing
of regular functions $f_1,\ldots,f_r$ on the affine scheme $Y$; we can
assume
that the linear span of the functions $f_1,\ldots,f_r$ is preserved by the
action
of $GL(n)$, defining a representation $GL(n)\arrow GL(r)$ over $k$. Then we
have
a $GL(n)$-equivariant morphism $\alpha:W\arrow A^r_k-0$
defined by
$$w\mapsto (f_1(w),\ldots,f_r(w)).$$
The subsets $\{ f_i\neq 0\} $ of $W$ are affine, and so $\alpha$ is an
affine morphism.

Define the affine group $\Aff_{r-1}$ as the semidirect product
 $(G_a)^{r-1}\rtimes GL(r-1)$. Then
we can identify $A^r-0$ with the homogeneous space $GL(r)/\Aff_{r-1}$.
Define $A$ as the pullback scheme:
$$\begin{CD}
A @>>> GL(r)\\
@VVV @VVV \\
W @>>> A^r-0
\end{CD}$$
Since $W$ is affine over $A^r-0$, $A$ is affine over $GL(r)$, and hence $A$
is an affine scheme. The group $GL(n)\times \Aff_{r-1}$ acts on $W$ and
$A^r-0$,
with $\Aff_{r-1}$ acting trivially, and it acts on $GL(r)$ via left
multiplication
by $GL(n)$ (using the representation $GL(n)\arrow GL(r)$) and right
multiplication
by $\Aff_{r-1}$. It follows that the pullback scheme $A$ also has an action
of $GL(n)\times \Aff_{r-1}$. We see from the diagram that we have an
isomorphism
of quotient stacks:
$$A/(GL(n)\times \Aff_{r-1})\cong W/GL(n).$$
So the given stack $X=W/GL(n)$ is the quotient of the affine scheme $A$
by the affine group scheme $GL(n)\times \Aff_{r-1}$. Thus (2) implies (3)
in Theorem \ref{main}.

To prove that (3) implies (2), we need the following well-known fact.

\begin{lemma}
\label{quot}
Every affine group scheme $G$ of finite type
over a field $k$ has a faithful representation
$G\arrow GL(n)$ such that $GL(n)/G$ is a quasi-affine scheme.
\end{lemma}

{\bf Proof. }To begin, let $G\arrow GL(n)$ be any faithful representation of
$G$.
By Chow, the homogeneous space $GL(n)/G$ is always
a quasi-projective scheme over $k$. More precisely,
there is a representation $V$ of $GL(n)$ and a $k$-point $x$ in the
projective
space $P(V^*)$ of lines in $V$ whose stabilizer is $G$ (\cite{DG}, p.~483).

Let the group $GL(n)\times G_m$ act on $V$ by the given representation of
$GL(n)$ and
by the action of $G_m$ by scalar multiplication. Then the stabilizer in
$GL(n)\times G_m$
of a point $y\in V$ lifting $x$ is isomorphic to $G$. That is, $(GL(n)\times
G_m)/G$ is a quasi-affine
scheme. Using the obvious inclusion $GL(n)\times G_m\arrow GL(n+1)$, the
quotient $GL(n+1)/(GL(n)\times G_m)$ is an affine scheme, the ``Stiefel
manifold'' of $(n+1)\times (n+1)$ matrices which are projections of rank 1.
So $GL(n+1)/G$ maps to the
affine scheme $GL(n+1)/(GL(n)\times G_m)$ with fibers the quasi-affine scheme
$(GL(n)\times G_m)/G$. More precisely, the fiber embeds
$(GL(n)\times G_m)$-equivariantly as a subscheme of the representation $V$,
by its construction.
 It follows
that $GL(n+1)/G$ is a quasi-affine scheme. \qed

We can now prove that (3) implies (2) in Theorem \ref{main}.
Let $X$ be the quotient stack
of an affine scheme $A$ of finite type over a field $k$
by the action of an affine group scheme $G$ of finite type over $k$. By
Lemma
\ref{quot},  there is a faithful representation $G\arrow GL(n)$ with
$GL(n)/G$ a quasi-affine scheme. Let $W$ be the $GL(n)$-bundle over
$X$ associated to the $G$-bundle $A\arrow X$, $W=(A\times GL(n))/G$.
 Then $W$ is an $A$-bundle over the quasi-affine scheme $GL(n)/G$. Since
$A$ is affine, it follows that $W$ is a quasi-affine scheme.

\section{From a stack to an algebraic space}
\label{stackspace}

We now begin the proof of the main part of Theorem \ref{main}, that (1)
implies (2). Thus, let $X$ be a normal 
noetherian stack whose stabilizer groups at closed points of $X$
are affine. Suppose that $X$
has the resolution property. We will show
that $X$ is isomorphic to the quotient stack $W/GL(n)$
for some quasi-affine scheme $W$ with an action of $GL(n)$. (In
Theorem \ref{main}, we only need $X$ to be normal for this part,
the proof that (1) implies (2). Normality
will be used in section \ref{scheme}.)
The following is a first step.

\begin{lemma}
\label{space}
Let $X$ be a noetherian stack (over $\Z$) whose stabilizer
groups at closed points of $X$ are affine.
Suppose that $X$ has the
resolution property. Then $X$ is isomorphic to the quotient stack
of some algebraic space $Z_1$ over $\Z$
by an action of the group $GL(n_1)$ for some $n_1$.
\end{lemma}

{\bf Proof. }Equivalently,
 we have to find a vector bundle $E_1$ on $X$ such that the total space $Z_1$
of the corresponding $GL(n_1)$-bundle over $X$ is an algebraic space.
As Edidin, Hassett, Kresch, and Vistoli observed, it is equivalent to require
that at every geometric point $x$ of $X$, the
action of the stabilizer group $G_x$ of $X$ on the fiber $(E_1)_x$ is
faithful
(\cite{EHKV}, Lemma 2.12).

By definition of a noetherian
stack, there is a smooth surjective morphism
from a noetherian
affine scheme $U$ to $X$ \cite{LMB}. We can think of the stack
$X$ as the quotient $U/R$
of $U$ by the groupoid $R:=U\times_{X}U$. In these terms, the defining
properties of a noetherian stack over $\Z$ are that $R$ is a separated
algebraic space over $\Z$ and that both projections $R\arrow U$
are smooth morphisms of finite type (\cite{LMB}, Proposition 4.3.1).
Let $G_U\arrow U$ be the stabilizer
group of $X$, that is, $G_U=R\times_{U\times U} U$. Here $G_U$
is a group in the category of algebraic spaces over $U$;
it is pulled back from the stabilizer group
$G:=X\times_{X\times X} X$ over $X$.  Since $G_U$ is a closed subspace
of $R$, it is separated of finite type over $U$.

A point of a stack $X$ is defined
in such a way that a point of $X$ is equivalent to an $R$-orbit in 
the underlying topological space of the scheme $U$. Another way to think
of a point is as an isomorphism class of substacks $\G $ of $X$ such that
$\G $ is a gerbe over some field $k$ (\cite{LMB}, Corollary 11.4);
explicitly, $\G $ is the quotient of the corresponding $R$-orbit
by the restriction of the groupoid $R$. To say that
$\G$ is a gerbe means that there is a field extension $F$ over $k$
such that $\G \times_k F$ is isomorphic to the classifying stack
of some group over $F$.
The set $|X|$ of points
of $X$ is given the quotient topology from $U$; in particular,
a closed point of $X$ can be identified with
a closed $R$-orbit in $U$ (\cite{LMB}, Corollary 5.6.1).

For any vector bundle $E$ on the stack $X$, the kernel of the $G$-action on
$E$
is a closed subgroup $H\subset G$ over $X$, which pulls
back to a closed subgroup $H_U\subset G_U$ over $U$. Given a finite sequence
of vector bundles $E_1,\ldots,E_n$ on $X$ with kernel subgroups
$H_1,\ldots,H_n\subset G$, the kernel subgroup of the direct sum
$E_1\oplus\cdots\oplus E_n$ is the intersection $H_1
\cap\cdots \cap H_n$. In this way,
we can repeatedly cut down the kernel subgroup by finding one vector
bundle on $X$ after another, and Lemma \ref{space} is proved if this
subgroup
eventually becomes the trivial group over $X$.

\begin{lemma}
\label{kernel}
Let $X$ be a noetherian stack (over $\Z$) which satisfies the
resolution property.
Let $x$ be a point of $X$ such that the stabilizer group $G$ of $X$
is affine at $x$. Then there is a vector
bundle $E$ on $X$ whose kernel subgroup is trivial at $x$.
\end{lemma}

{\bf Proof. }As explained above, $x$ corresponds to a substack $\G $
of $X$ which is a gerbe over some field $k$. Since $X$ is locally
noetherian, there is a finite extension $F$ of $k$ such that
$\G\times_k F$ is isomorphic to 
the classifying stack of a group $G_s$ over $F$, by \cite{LMB},
11.2.1 and Theorem 11.3. The assumption means that $G_s$ is affine over $F$.
Moreover, $G_s$ is of finite type over $F$, since $G\arrow X$ is of
finite type.
Therefore $G_s$ has a faithful representation over $F$. We can view
such a representation
as a vector bundle on the gerbe $\G \times_k F$. Its direct
image to $\G $ is a vector bundle $C_0$ on $\G $ whose pullback 
to $\G \times_k F$ is a faithful representation of the group $G_s$.
So the kernel subgroup of $C_0$ over $\G $ is trivial.

Let $i:\G \arrow X$ denote the inclusion. The direct image $i_*C_0$
is a quasi-coherent sheaf on $X$, and therefore a direct
limit of coherent sheaves on $X$ (\cite{LMB}, Proposition 13.2.6 and
Proposition 15.4). Since $X$ has the resolution property, each
of these coherent sheaves is a quotient of a vector bundle on $X$.
Since $C_0=i^*i_*C_0$ and $C_0$ is coherent, one of these vector bundles
$E$ on $X$ must restrict to a vector bundle on $\G $ which maps onto $C_0$.
It follows that the kernel subgroup of $E$ over $\G $ is trivial. \qed

We return to the proof of Lemma \ref{space}. We are given that the
stabilizer group of $X$ at each closed point of $X$ is affine. By Lemma
\ref{kernel}, it follows that for each closed point $x$ of $X$, there
is a vector bundle $E$ on $X$ whose kernel subgroup $H\arrow X$
is trivial at $x$.
This does not imply that the kernel subgroup of $E$ is trivial in
a neighborhood of $E$. Nonetheless,
since the morphism $H\arrow X$ has finite type, the dimensions of fibers
make sense and are upper semicontinuous. (Indeed, it suffices
to check this for the pulled-back group $H_U\arrow U$, and to consider
an etale covering of the algebraic space $H_U$ by a scheme; then
we can refer to EGA IV.13.1.3 \cite{EGA}.)
Therefore the group $H\arrow X$ is quasi-finite (that is, of
finite type and with
finite fibers) over some neighborhood of the point $x$. The space
$|X|$ of points of $X$ is a ``sober'' noetherian topological space
(every irreducible closed subset of $|X|$ has a unique generic point)
by \cite{LMB}, Corollary 5.7.2. It follows that 
every open subset of $|X|$ which contains
all the closed points must be the whole space. Since
$|X|$ is also quasi-compact,
there are finitely many vector bundles
$E_1,\ldots,E_n$ on $X$ such that the direct sum $E_1\oplus
\cdots \oplus E_n$ has kernel group which is quasi-finite over all
of $X$. Let $E$ now denote this direct sum.
In particular, the kernel subgroup $H\arrow X$ of $E$ is
affine over every point of $X$.

Since the kernel group $H_U\arrow U$ is quasi-finite, it
is finite over a dense open subset $V$ of $U$
(\cite{Hartshorne}, exercise II.3.7). By Lemma \ref{kernel}, for every
point $x\in V$, there is a vector bundle $F$ on $X$ whose
kernel subgroup is trivial at $x$. Since $H_{E\oplus F}$ is a closed
subgroup scheme of $H_E$, it is finite over $V$, while also being trivial
at $x$. Therefore $H_{E\oplus F}$ is trivial over some neighborhood of
$x$.
By quasi-compactness of $V$, there is a vector bundle on $X$
(again to be called $E$) whose kernel subgroup is
quasi-finite
over $U$ and trivial over $V$. Then this kernel subgroup
will be finite over a larger open
subset of $U$, containing a dense open subset of $U-V$, and so we can repeat
the process. By noetherian induction, we end up with a vector bundle on
$X=U/R$ with trivial kernel subgroup over all of $U$. \qed

\section{From an algebraic space to a scheme}
\label{scheme}

In this section, we will complete the proof of Theorem \ref{main}.

Let $X$ be a normal
noetherian stack $X$ with affine stabilizer groups at closed points of $X$.
Suppose that $X$ satisfies the resolution property.
We will show that $X=W/GL(n)$ for some
quasi-affine scheme $W$ and some $n$. By Lemma \ref{space}, we know that
there is
a $GL(n_1)$-bundle $Z_1$ over $X$ for some $n_1$ which is at least an
algebraic space. Since $X$ is normal, $Z_1$ is normal. We use
the name $E_1$ for the vector bundle on $X$ that corresponds to
the $GL(n_1)$-bundle $Z_1$.

By Artin, since $Z_1$ is a normal noetherian algebraic
space, it is the coarse geometric quotient of some normal scheme $A$
by the action of a finite group $G$ (\cite{KollarJDG}, 2.8;
\cite{LMB}, 16.6.2).
Moreover, the morphism $\pi:A\arrow Z_1$ is finite. (In general,
the quotient morphism even of an affine noetherian scheme by a finite
group need not be a finite morphism, by Nagata \cite{Nagata}, but the
situation
here is better because we know $Z_1$ is noetherian to start with.)
Let $\pi:A\arrow Z_1$ be the corresponding morphism.
Let $U_1,\ldots,U_r$ be an open affine covering
of the scheme $A$. Let $S_i$ be the closed subset $A-U_i$, which we give
the reduced subscheme structure. Let $I_{S_i}$ be the corresponding ideal
sheaf
(the kernel of $O_A\arrow O_{S_i}$), and let $C$ be the coherent sheaf
$C=\oplus_{i=1}^r I_{S_i}$ on $A$. Then $D:=\pi_*C$ is a coherent sheaf on
$Z_1$
because $\pi:A\arrow Z$ is proper.

In order to use the resolution property for $X$ again, we need a suitable
coherent sheaf on $X=Z_1/GL(n)$. That will be supplied by the following
lemma.

\begin{lemma}
\label{equi}
Let $Z_1$ be a noetherian algebraic space (over $\Z$). 
Let $G$ be a flat affine group scheme over $\Z$ or over a field which
acts on $Z_1$.
Then any coherent sheaf on $Z_1$ is a quotient of some $G$-equivariant
coherent sheaf on $Z_1$.
\end{lemma}

{\bf Proof. }The morphism $\alpha$ from $Z_1$ to the quotient stack $Z_1/G$
is affine.
So, for every coherent sheaf $D$ on $Z_1$, the natural map
$$\alpha^*\alpha_*D\arrow D$$
is surjective.
Here $\alpha_*D$ is a quasi-coherent sheaf on $Z_1/G$. Thus, we have
exhibited
the coherent sheaf $D$ as the quotient of a $G$-equivariant quasi-coherent
sheaf
on $Z_1$. By Laumon and Moret-Bailly, every quasi-coherent sheaf on a
noetherian
stack is the filtered direct limit of its coherent subsheaves (\cite{LMB},
Proposition
15.4). So $D$ is in fact the quotient of some $G$-equivariant coherent sheaf
on $Z_1$.
\qed

Before continuing with the proof of Theorem \ref{main}, note the following
corollary which could be useful in checking the resolution property
in examples. For example, to prove the resolution property for
all coherent sheaves on a toric variety, it suffices to prove it
for the equivariant coherent sheaves. Klyachko's algebraic
description of the equivariant vector bundles on a toric variety
should be useful for the latter problem \cite{Klyachko}.

\begin{corollary}
\label{quotres}
Let $Z_1$ be a noetherian algebraic space (over $\Z$).
Let $G$ be a flat affine group scheme of finite type
over $\Z$ or over a field which
acts on $Z_1$. If the resolution property holds for the stack $Z_1/G$,
then it holds for $Z_1$.
\end{corollary}

{\bf Proof of Corollary \ref{quotres}. }Every coherent sheaf
on $Z_1$ is a quotient of a $G$-equivariant coherent sheaf by
Lemma \ref{equi}, which in turn is a quotient of a $G$-equivariant
vector bundle on $Z_1$ by the resolution property for $Z_1/G$. \qed

We now return to the proof of Theorem \ref{main}.
We apply Lemma \ref{equi} to the algebraic space $Z_1$ with $X=Z_1/GL(n_1)$
and the coherent sheaf $D=\pi_*C$ on $Z_1$ defined above. It follows that
$\pi_*C$
is the quotient of some $GL(n_1)$-equivariant coherent sheaf $D_1$ on $Z_1$.
 Since the stack $X$ has
the resolution property, there is a vector bundle $E_2$ on $X$ which maps
onto $D_1$, viewed as a coherent sheaf on $X$. Thus, writing $\alpha$
for the $GL(n_1)$-bundle $Z_1\arrow X$, we have surjections
$$\alpha^*E_2\surj D_1\surj \pi_*C$$
on $Z$. Since the morphism $\pi:A\arrow Z_1$ is finite,
it is affine. So the natural
map $\pi^*\pi_*C\arrow C$ of sheaves 
over $A$ is surjective. Thus, we have found
a vector bundle $E_2$ on $X$ whose pullback $\pi^*\alpha^*E_2$
to the scheme $A$ maps onto the coherent sheaf $C=\oplus_{i=1}^nI_{S_i}$.

Let $Z_2$ be the $GL(n_2)$-bundle over $X$ associated to the vector bundle
$E_2$. Define $W$ and $Y$ as the indicated pullbacks:
$$ \begin{CD}
W @>{\beta}>{GL(n_2)}>  A\\
@VVV  @VV{\pi}V\\
Y @>>{GL(n_2)}>  Z_1\\
@VVV  @V{\alpha}V{GL(n_1)}V\\
Z_2 @>>{GL(n_2)}>  X
\end{CD} $$

Thus $W$ is the $GL(n_2)$-bundle over the scheme $A$ associated to
the vector bundle $\pi^*\alpha^*E_2$, which we know maps onto the coherent
sheaf $C=\oplus_{i=1}^n I_{S_i}$. It follows that the scheme $W$ is
quasi-affine,
by the following argument. By construction of $W$, the vector bundle
$\pi^*\alpha^*E_2$ on $A$ pulls
back to the trivial bundle on $W$. Let $\beta:W\arrow A$ denote
this $GL(n_2)$-bundle. Then
the coherent sheaf $\beta^*(\oplus_{i=1}^n
I_{S_i})=\oplus_{i=1}^nI_{\beta^{-1}(S_i)}$ on $W$ is spanned by its global
sections. This means that for each $1\leq i\leq n$, the open subset
$W-\beta^{-1}(S_i)$ is the union of its open subsets of the form
$\{ f_{ij}\neq 0\} $, for certain regular functions $f_{i1},\ldots,f_{ir}$
on $W$
which vanish on $\beta^{-1}(S_i)$. Moreover, each subset $S_i$ was chosen so
that
$A-S_i$ is an affine scheme. Since $\beta:W\arrow A$ is an affine morphism,
$W-\beta^{-1}(S_i)$ is also an affine scheme. It follows that the open subsets
$\{ f_{ij}\neq 0\}$ of $W-\beta^{-1}(S_i)$
are affine. Thus, using all $i$ and $j$,
$W$ is the union of open affine subschemes of the form $\{ f\neq 0\}$ for
regular functions $f$ on $W$. This means that the scheme $W$ is quasi-affine
(\cite{EGA}, Proposition II.5.1.2).

By the pullback diagram above, $Y$ is a $GL(n_1)\times GL(n_2)$-bundle over
$X$.
Because $Y$ is a $GL(n_2)$-bundle over the algebraic space $Z_1$, $Y$
is an algebraic space. Moreover, as the pullback diagram shows, we have
a finite surjective morphism from the quasi-affine scheme $W$ to $Y$.
In general, this does not imply that $Y$ is a quasi-affine scheme,
as Grothendieck observed (\cite{EGA}, Remark II.6.6.13); for example,
there is a non-quasi-affine scheme whose normalization is quasi-affine.

We know, however, that $Y$ is normal, and that $Y$ is the coarse geometric
quotient of the quasi-affine scheme $W$ by a finite group $G$. It follows
by the usual construction of quotients by finite group actions
that $Y$ is a quasi-affine scheme. Namely, since $W$ is
quasi-affine, every finite subset of $W$ is contained in an affine open
subset of the form $\{ f\neq 0\}$ for some regular function $f$ on $W$.
 In particular, each $G$-orbit in $W$
is contained in an affine open subset. This subset can be taken to be
of the form $\{ f\neq 0\}$ for some $G$-invariant function $f$,
by taking the product of the translates of a given function on $W$.
Then we can define the geometric quotient $Y$ of $W$ by $G$
as a scheme, the union of open subsets $\Spec O(U)^G$ corresponding
to these affine open subsets $U$ of $W$; a reference that works in
this
generality is SGA 3 (\cite{Gabriel},
Theorem 4.1). Finally, the scheme $Y$ thus defined
is quasi-affine because it is an open subset
of the affine scheme $\Spec O(W)^G$. (The rings $O(W)$ and $O(W)^G$ may not be
noetherian, but that does not matter for the purpose of proving
that $Y$ is quasi-affine.)

We can then consider the $GL(n_1+n_2)$-bundle over $X$ associated
to the vector bundle $E_1\oplus E_2$. Its total space is a bundle over
$Y$ with fiber the affine scheme $GL(n_1+n_2)/(GL(n_1)\times GL(n_2))$
(a Stiefel manifold as in section \ref{easy}),
and hence is a quasi-affine scheme. Theorem \ref{main} is proved. \qed

\section{Proof of Proposition \ref{affine}}

We now prove Proposition \ref{affine}. That is, let $X$ be a noetherian stack
(over $\Z$) with affine stabilizer groups at closed points of $X$.
Suppose that $X$ has the resolution
property. We will show that
the diagonal morphism $X\arrow X\times_{\Z}X$ is affine.
In this section, I will write $X\times Y$ to mean $X\times_{\Z}Y$.

Suppose first that $X$ is normal. By Theorem \ref{main}, $X$ is isomorphic
to the quotient of some quasi-affine scheme $W$ by an action of $GL(n)$.
A quasi-affine scheme $W$ is separated (that is, the diagonal morphism
$W\arrow W\times W$ is a closed embedding). In particular, $W$ has
affine
diagonal. Since $GL(n)$ is an affine group scheme over $\Z$, it follows that
$X=W/GL(n)$ has affine diagonal. This is an easy formal argument, as follows.
For brevity, let us write $G$ for $GL(n)$. Since $W\arrow W/G$
is faithfully flat, to show that $W/G\arrow W/G\times W/G$ is
affine is equivalent to showing that the pulled-back map over $W\times W$
is affine, that is, that $G\times W\arrow W\times W$,
$(g,w)\mapsto (w,gw)$, is affine. But that map is the composition
of the map $G\times W \arrow G\times W\times W$
by $(g,w)\mapsto (g,w,gw)$, which is a pullback of the diagonal map
of $W$ and hence is affine, with the projection map
$G\times W\times W
\arrow W\times  W$ which is affine since $G$ is affine.

For an arbitrary noetherian stack $X$ with affine stabilizer groups
at closed points,
the proof of Theorem \ref{main} works until the last step. We find
that $X$ is isomorphic to $Y/GL(n)$ for some noetherian algebraic space $Y$
which admits a finite surjective morphism $W\arrow Y$ from a
quasi-affine scheme $W$. The point now is that $Y$ is separated since $W$ is,
by the following lemma. In particular, $Y$ has affine diagonal, and
so the stack $X=Y/GL(n)$ has affine diagonal. Proposition \ref{affine}
is proved. \qed

\begin{lemma}
Let $f:X\arrow Y$ be a proper surjective morphism of noetherian algebraic
spaces.
If $X$ is separated, then $Y$ is separated.
\end{lemma}

{\bf Proof. }We have the commutative diagram
$$\begin{CD}
X @>>> X\times X \\
@VfVV @VVV \\
Y @>g>> Y\times Y
\end{CD}$$
The map $X\arrow X\times X$ is proper since $X$ is separated,
and $X\times X\arrow Y\times Y$ is proper, so the composition
$X\arrow Y\times Y$ is proper. By the defining properties
of a noetherian stack over $\Z$ (as in section \ref{stackspace} or
\cite{LMB}),
the diagonal morphism $g:Y\arrow Y\times Y$ is separated and of finite
type. We would like to conclude
that $g:Y\arrow Y\times Y$ is proper, that is,
$Y$ is separated, by EGA II.5.4.3 \cite{EGA}: ``Let $f:A\arrow B$,
$g:B\arrow C$ be morphisms of schemes such that $g\circ f$ is
proper. If $g$ is separated of finite type and $f$ is surjective,
then $g$ is proper.'' We have algebraic spaces rather than schemes here,
but the same proof works, as follows. By definition of properness for a map
of algebraic spaces (\cite{Knutson}, Definition II.7.1),
since we know that $g$ is separated of finite type,
we only have to show that $g$ is universally closed. The hypotheses
pull back under arbitrary morphisms of algebraic spaces, so it suffices
to show that the morphism $g$ is closed, that is, that the images
of closed sets are closed. But this is clear from
surjectivity of $f$ together with the fact that the morphism
$g\circ f$ is closed. \qed

\section{Varieties with quotient singularities: proof of Theorem \ref{orb2}}

We now prove Theorem \ref{orb2}. Thus,
let $X$ be a smooth Deligne-Mumford stack over a field $k$. Suppose that
$X$ has finite stabilizer group and that the stabilizer group is generically
trivial. Let $B$ be the Keel-Mori coarse moduli space of $X$ \cite{KM}.
Finally, suppose that
the algebraic space
$B$ is a scheme with affine diagonal (for example a separated scheme).
We will show that the stack 
$X$ has the resolution property.

We need the following important property of the Keel-Mori space,
which I have stated in the full generality in which that space is defined.

\begin{lemma}
\label{km}
Let $X$ be a stack of finite type over a locally noetherian base scheme
$S$. Suppose that $X$ has finite stabilizer group, so that there
is a Keel-Mori quotient space $B$. Then the map $X\arrow B$ is proper.
\end{lemma}

This is a more general version of Keel and Mori \cite{KM}, 6.4,
which in turn is modeled on Koll\'ar \cite{Kollarquot}, 2.9.

{\bf Proof. }Here properness for a map of stacks is defined in
\cite{LMB}, Chapter 7. In the case at hand, there is a finite
surjective morphism from a scheme to $X$ (\cite{EHKV}, Theorem 2.7).
As a result, there is a valuative criterion for properness
using only discrete valuation rings (\cite{LMB}, Proposition 7.12),
in which we ask
for a lift after a suitable ramified extension of DVRs.

The problem is local over $S$, so we can assume that $S$ is noetherian.
So $X$ is noetherian. Therefore
we can find a smooth surjective morphism from an affine scheme $U$
to $X$. As in section \ref{stackspace}, $X$ is the quotient stack
of $U$ by the groupoid $R:=U\times_X U$. By property 1.8US of the Keel-Mori
quotient, the map $U\arrow B$ is a universal submersion
(as defined in EGA IV.15.7.8 \cite{EGA}). Therefore, as Koll\'ar
explains, if $T$ is the spectrum of a DVR and $u:T\arrow B$
is a morphism then there is a dominant morphism, which we can assume
to be finite, from another DVR $T'$ to $T$, such that the composition
$T'\arrow T\arrow B$ lifts to a map $\overline{u}: T'\arrow U$.
In the valuative criterion for properness of $X\arrow B$, we are also given
a lift $v$ of $u$ restricted to the general point $t_g$ of $T$,
$v:t_g\arrow U$. By property 1.8G of the Keel-Mori quotient,
$U(\xi)/R(\xi)\arrow B(\xi)$ is a bijection for every geometric
point $\xi$, and so $\overline{u}|_{t_g'}$ and $v$ become
equivalent under the groupoid $R$ after base extension to another
DVR $T''$ which is finite over $T'$. Since $\overline{u}$ is defined
on all of $T''$, this checks the valuative criterion: the morphism
$X\arrow B$ is proper. \qed

We return to the proof of Theorem \ref{orb2}. 
Since $X$ is a smooth Deligne-Mumford stack with trivial
generic stabilizer, it is the quotient of some algebraic space $Z_1$
by an action of $GL(n_1)$ over $k$, by Edidin, Hassett, Kresch,
and Vistoli (\cite{EHKV}, Theorem 2.18). (In characteristic zero,
this is essentially Satake's classical observation that an orbifold
with trivial generic stabilizer is a quotient of a manifold by a
compact Lie group,
using the frame bundle corresponding to the tangent bundle
(\cite{Satake}, p.~475).)

By Laumon and 
Moret-Bailly (generalized by Edidin, Hassett, Kresch, and Vistoli,
as used in the above proof), there is a finite surjective morphism
from a scheme $A$ to the Deligne-Mumford stack $X$ (\cite{LMB},
Theorem 16.6). Some special cases of this result were known before,
by Seshadri \cite{Seshadri}, 6.1, and Vistoli \cite{Vistoli}, 2.6.
By Lemma \ref{km}, the morphism $X\arrow B$ is proper, and
so the composition $A\arrow B$ is proper. It is clearly also
a quasi-finite morphism of algebraic spaces, and so it is finite,
in particular affine.
Define $Z_A$
by the following pullback diagram.
$$ \begin{CD}
Z_A @>{}>{}>  Z_1\\
@V{}VV  @VV{GL(n_1)}V\\
A @>>{}>  X\\
@.  @V{}V{}V\\
@.   B
\end{CD} $$
Since $Z_A\arrow A$ is a $GL(n_1)$-bundle and $A\arrow B$ is finite,
both morphisms are affine, and so the composition $Z_A\arrow B$ is
affine. Since $Z_A\arrow Z_1$ is a finite surjective morphism
of algebraic spaces over $B$, it follows from Chevalley's theorem
for algebraic spaces (\cite{Knutson}, III.4.1) that the morphism $Z_1\arrow B$
is affine.

Since $B$ is a scheme and the morphism $Z_1\arrow B$ is affine,
the smooth algebraic space 
$Z_1$ is a scheme. Likewise, since $B$ has affine diagonal
and $Z_1\arrow B$ is affine, $Z_1$ has affine diagonal. To see this,
write the diagonal morphism $Z_1\arrow Z_1\times_k Z_1$ as the composition
of two maps. The first is
$Z_1\arrow Z_1 \times_B Z_1$, which is a closed embedding and hence
affine since $Z_1\arrow B$ is affine and hence separated. Next is
$Z_1\times_B Z_1\arrow Z_1\times_k Z_1$, which is a pullback of the affine
morphism $B\arrow B\times_k B$ and hence is affine.

Thus, $Z_1$ is a smooth scheme with affine diagonal. So $Z_1$ has
an ample family of line bundles, by
Brenner and Schr\"oer (or Kleiman and Illusie, in the separated case),
as mentioned in section \ref{history}.
Then, by Theorem \ref{thom}, due to Thomason, the
quotient stack $X = Z_1/GL(n)$ has the resolution property. \qed

\section{The resolution property and $K$-theory}
\label{counterexamples}

As mentioned in section \ref{history}, if a stack $X$ has the
resolution property, then the natural map $K_0^{\naive}X\arrow K_0X$
is an isomorphism. In particular, every perfect complex on $X$ is equivalent
in the Grothendieck group $K_0X$ of perfect complexes
to a difference of vector bundles. We will present two examples to show
that the converses to these statements are false.
First, we state a positive result for smooth
schemes.

\begin{proposition}
Let $X$ be a smooth scheme of finite type over a field. The following
are equivalent.

(1) $X$ has affine diagonal. For a scheme, as here, it is equivalent to say
that the intersection of any
two affine open subsets of $X$ is affine.

(2) $X$ is a scheme with an ample family of line bundles.

(3) $X$ has the resolution property.

(4) The natural map from $K_0^{\naive}X$ to $K_0X=G_0X$ is surjective.
Equivalently, every coherent sheaf on $X$
is equivalent in the Grothendieck group $G_0X$ to a difference of
vector bundles.
\end{proposition}

This is fairly easy, but perhaps suggestive.
From my point of view, the interesting equivalence here is between (1)
and (3), because one would often like to know whether the resolution
property holds, and it is usually easy to check whether a scheme or stack
has affine diagonal. We can hope that the equivalence between (1) and (3)
holds in much greater generality; see Questions 1 and 2 in the introduction.
In more general situations, property (2) will imply the resolution property
but will definitely not be equivalent to it, as can be seen from several
examples in section \ref{history}. Finally, the end of this section
will present two examples showing that property (4) does not imply
the resolution property in more general situations, for example for
smooth algebraic spaces.

{\bf Proof. }The implications $(1)\imp (2) \imp (3) \imp (4)$ are
discussed in section \ref{history}. First, Brenner and Schr\"oer
observed that (1) implies (2), that is, that a smooth scheme with affine
diagonal has an ample family of line bundles. The proof is the same
as in the case of a smooth separated scheme, due to Borelli
and Illusie. Next,
Kleiman and Illusie proved that (2) implies (3).
Finally, Thomason proved that $K_*^{\naive}X\arrow K_*X$
is an isomorphism when $X$ has the resolution property. The special
case that (3) implies (4) is particularly simple: using the resolution
property, every coherent sheaf has a resolution by vector bundles,
which can be stopped after finitely many steps because the scheme $X$
is regular.

It remains to show that (4) implies (1). That is, if $X$ is a smooth
scheme that does not have affine diagonal, we will find a coherent sheaf $C$
on $X$ whose class in $K_0X=G_0X$ is not a difference of vector
bundles. The following proof extends an argument by Schr\"oer
and Vezzosi (\cite{SV}, Proposition 4.2).

The assumption that $X$ does not have affine diagonal means that $X$
has open affine subsets $U$ and $V$ such that $U\cap V$ is not affine.
As is well known, the complement of an irreducible
divisor $D$ in a smooth affine variety $U$
is affine. Indeed, $D$ is Cartier because $U$ is smooth, and so the
inclusion from $U-D$ into $U$ is an affine morphism, which implies that
$U-D$ is affine since $U$ is. As a result, if $U-(U\cap V)$ contains
any irreducible divisor, we can remove it from $U$ without
changing the properties we have stated for $U$, and likewise for $V$.
Thus we can assume that $U-(U\cap V)$ and $V-(U\cap V)$ have codimension
at least 2.

Since $U$ and $V$ are smooth, in particular normal, it follows
that the restriction maps from $O(U)$ or $O(V)$ to $O(U\cap V)$
are isomorphisms. Since $U$ and $V$ are affine, this means that both
$U$ and $V$ are isomorphic to $\Spec O(U\cap V)$. Let $S$ be the closure
in $X$ of $U-(U\cap V)$, and let $T$ be the closure in $X$ of
$V-(U\cap V)$. We give these closed subsets the reduced scheme
structure. Suppose that
the coherent sheaf $O_S$ on $X$ is a difference of vector bundles
$E-F$ in the Grothendieck group $G_0X$; we will derive a contradiction.

After shrinking $U$ to a smaller affine open neighborhood of the
generic point of an irreducible component of $S$,
and shrinking $V$ to the corresponding neighborhood of
the generic point of a component of $T$, we can
assume that the vector bundles $E$ and $F$ are trivial on both $U$ and $V$.
So each of these vector bundles is described up to isomorphism on $U\cup V$
by an attaching map $U\cap V\arrow GL(n)$. Since $U-(U\cap V)$
has codimension at least 2, every such map extends to $U$. It follows 
that $E$ and $F$ are in fact trivial on $U\cup V$. Thus, our
assumption implies that the class of $O_S$ in $G_0X$ restricts to zero
in $G_0(U\cup V)$.

Let $S_U=S\cap U=U-(U\cap V)$ and $T_V=T\cap V=V-(U\cap V)$;
these are disjoint nonempty 
closed subsets of $U\cup V$, and they are isomorphic.
Consider the localization sequences in $G$-theory, due to Quillen
\cite{Quillen}:
$$ \begin{CD}
G_1(U\cap V)@>>> G_0S_U@>>> G_0U@>>> G_0(U\cap V)@>>> 0\\
G_1(U\cap V)@>>> G_0T_V@>>> G_0V@>>> G_0(U\cap V)@>>> 0\\
G_1(U\cap V)@>>> G_0(S_U{\textstyle \, \coprod \,} T_V)
@>>> G_0(U\cup V)@>>> G_0(U\cap V)@>>> 0
\end{CD} $$
Since the inclusions of $U\cap V$ into $U$ and into $V$ are isomorphic,
any element of $G_1(U\cap V)$ has the same image in $G_0S_U$ as in
$G_0T_V$, with respect to the isomorphism of $S_U$ with $T_V$. Furthermore,
the class of $O_S$ in $G_0S_U$ is not zero. So the class of $O_S$ in
$G_0(S_U\coprod T_V)=G_0S_U\oplus G_0T_V$ is not in the image of
$G_1(U\cap V)$. Thus the class of $O_S$
in $G_0(U\cup V)$ is not zero, contradicting the previous paragraph.
So in fact the class of the coherent sheaf $O_S$ in $G_0X$
is not a difference of vector bundles. We have proved that (4) implies (1).
\qed

{\bf Example 1. }There is a smooth algebraic space $Z$ such that
$K_0^{\naive}Z\arrow K_0Z$ is surjective, that is, every perfect complex on $Z$
is equivalent in $K_0Z$ to a difference of vector bundles,
but $Z$ does not have the resolution property.

In conformity with Question 1 in the introduction, the space
we construct will not have affine diagonal. Let $Y_r$ be
the smooth non-separated scheme $Y=A^r\cup_{A^r-0}A^r$,
$r\geq 2$, over some field $k$ of characteristic not 2.
The algebraic space $Z_r$ will be the
quotient of $Y_r$
by a free action of the group $\Z/2$,
acting by $-1$ on $A^r$ and
switching the two origins. The algebraic space $Z_r$
is a well known example, described and illustrated by Artin \cite{Artin}
and named by Koll\'ar
a bug-eyed cover \cite{Kollarbug}.
Its best known property is that it is not locally
separated at the image of the origin, and therefore not a scheme.

For $r=1$, the algebraic space $Z_1$ has the resolution property.
Indeed, $Y_1$ is a smooth scheme with affine diagonal,
and so it has an ample family of line bundles by the result of
Brenner and Schr\"oer mentioned in section \ref{history}. By
Theorem \ref{thom},
it follows that $Z_1=Y_1/(\Z/2)$ has the resolution property.
One gets a more direct proof by observing that
the scheme
$Y_1=A^1\cup_{A^1-0}A^1$ is the quotient of $A^2-0$ by
the diagonal torus $G_m \subset SL(2)$.
Likewise, the algebraic space $Z_1$ is the quotient of the quasi-affine scheme
$A^2-0$ by the normalizer of $G_m$ in $SL(2)$.
(This normalizer is a non-split extension of $\Z/2$ by $G_m$.) Then it is
immediate from Theorem \ref{thom} that $Y_1$ and $Z_1$ have the
resolution property.

We now consider $r\geq 2$. Since $Z_r$ does not have affine diagonal,
we know by Proposition \ref{affine} that $Z_r$ does not have the resolution
property, as we will  see more explicitly below.
To compute the group $K_0^{\naive}Z_r$, we need to describe
the vector bundles on $Z_r$. Let $\sigma:A^r-0\arrow A^r-0$ be
multiplication by $-1$. A vector bundle $E$ on $Z_r$ is a vector
bundle $E$ on $A^r$ together with an isomorphism
$$f:E\overset{\cong}{\arrow} \sigma^*E$$
over $A^r-0$, such that the composition
$$E\overset{f}{\arrow} \sigma^*E\overset{\sigma^*f}{\arrow}
\sigma^*\sigma^*E=E$$
is the identity on $A^r-0$. Since $r\geq 2$, $f$ extends uniquely to a map
$f:E\arrow \sigma^*E$ over all of $A^r$, and the above composition
is the identity over $A^r$ because this is true over $A^r-0$. Thus
the category of vector bundles over $Z_r$ is equivalent to that
of $\Z/2$-equivariant vector bundles over $A^r$, with $\Z/2$ acting
on $A^r$ by multiplication by $-1$. So
\begin{align*}
K_0^{\naive}Z_r&\cong K_0^{\naive}(A^r/\Z/2)\\
&\cong K_0(A^r/\Z/2)\\
&\cong \text{Rep}(\Z/2)\\
&\cong \Z^2.
\end{align*}
Here $A^r/\Z/2$ denotes the quotient stack of $A^r$ by $\Z/2$.
 Its naive $K$-theory coincides
with its true $K$-theory because it has the resolution property, by
Theorem \ref{thom}. The calculation that $K_0(A^r/\Z/2)$ is isomorphic
to the representation ring of $\Z/2$ follows from the
homotopy invariance of equivariant algebraic $K$-theory, also proved
by Thomason (\cite{Thomasongroup}, 4.1).

Next, we compute the true $K$-theory $K_0Z_r$, which is isomorphic
to the Grothendieck group $G_0Z_r$ of coherent sheaves because
the algebraic space $Z_r$ is smooth over $k$. By the previous paragraph,
we can identify the map $K_0^{\naive}Z_r\arrow K_0Z_r$ with the pullback map
$G_0(A^r/\Z/2)\arrow G_0Z_r$ associated to the obvious flat morphism
from $Z_r$ to the quotient stack $A^r/\Z/2$.
We have exact localization sequences, by Thomason's paper on equivariant
$K$-theory (\cite{Thomasongroup}, 2.7):
$$\begin{CD}
G_0(\text{point}/\Z/2)@>>> G_0(A^r/\Z/2)@>>> G_0((A^r-0)/\Z/2)@>>> 0.\\
@VVV @VVV @VVV\\
G_0((\text{2 points})/\Z/2)@>>> G_0Z_r@>>> G_0((A^r-0)/\Z/2)@>>> 0.\\
\end{CD}$$
The left vertical map sends $\text{Rep}(\Z/2)\cong \Z^2$ to
$G_0(\text{point})=\Z$ by the rank; in particular, it is surjective.
The right vertical map is an isomorphism, and so the center vertical
map is surjective. This means that $K_0^{\naive}Z_r\arrow
K_0Z_r$ is surjective, as we want.

We can also compute $K_0Z_r$ ($=G_0Z_r$) explicitly, using the above
diagram.  By a Koszul resolution, the pushforward map
$$\Rep(\Z/2)=G_0(\text{point}/\Z/2)\arrow G_0(A^r/\Z/2)=\Rep(\Z/2)$$
is multiplication by $\lambda_{-1}V:=\sum (-1)^i\Lambda^iV$, where $V$
denotes the representation of $\Z/2$ on $A^r$. Therefore $K_0Z_r$
is the quotient of $\Rep(\Z/2)$ by the relation
$$W\cdot \lambda_{-1}V=(\text{dim }W)\cdot \lambda_{-1}V$$
for all $W\in\Rep(\Z/2)$. We thereby compute
that $K_0Z_r$ is the quotient of
$K_0^{\naive}Z_r=\text{Rep}(\Z/2)
=\Z\oplus \Z u$, where $u$ is the nontrivial 1-dimensional representation
of $\Z/2$, by the relation $2^r(1-u)=0$, so that $K_0Z_r$ is isomorphic
to $\Z\oplus \Z/2^r$. In particular, we see again that
the resolution property fails for $Z_r$, because the map $K_0^{\naive}Z_r
\arrow K_0Z_r$ is not an isomorphism.
Explicitly, let $O_{A^r}$ and $L$ denote the $\Z/2$-equivariant
line bundles on $A^r$ associated to the trivial and the nontrivial
1-dimensional representations of $\Z/2$. Let $K$ be the coherent sheaf
$$K=\ker(O_{A^r}\oplus L\arrow O_0),$$
where both $O_{A^r}$ and $L$ map onto $O_0$. Then $K$ is not
a $\Z/2$-equivariant coherent sheaf on $A^r$, but it is $\Z/2$-equivariant
outside the origin, and so it corresponds to a coherent sheaf on $X_r$.
The sheaf $K$ is not 
a quotient of a vector bundle on $X_r$.

It is amusing to observe that,
over the complex numbers, the non-separated scheme $Y_r=A^r\cup_{A^r-0}
A^r$ is weak homotopy equivalent to the sphere $S^{2r}$, and the
quotient algebraic space $Z_r$ is weak homotopy equivalent to real
projective
space $\R\P^{2r}$. The true $K$-group $K_0Z_r=\Z\oplus \Z/2^r$
maps isomorphically
to the topological $K$-group  $K^0_{\text{top}}\R\P^{2r}$, as computed by
Atiyah (\cite{Atiyahbook}, p.~107). This suggests that the $K$-theory of any
stack
with affine stabilizer group
should be closely related to its topological or etale $K$-theory. The
relation will not always be visible on the level of $K_0$, but rather in
the groups $K_i$ with $i$ large. Precisely,
there should be an isomorphism
$$K_*(X;\Z/l^{\nu})[\beta^{-1}]\arrow K^*_{\text{et}}(X;\Z/l^{\nu}),$$
where $\beta$ denotes the Bott map. In fact, one problem here
is to define the groups on the right. For $X$ a locally separated algebraic
space of finite type over a reasonable base scheme,
or more generally for quotient stacks $X/G$ of such
an algebraic space by a linear algebraic group, Thomason proved analogous
results in $G$-theory (\cite{Thomasonequi}, 3.17 and Theorem 5.9),
which coincides with $K$-theory for regular stacks.
Some of his results apply to the above example $Z_r=Y_r/\Z/2$.

{\bf Example 2. }There is a smooth Deligne-Mumford stack $X$ for which
the natural map $K_0^{\naive}X\arrow K_0X$ is an isomorphism, but $X$
does not have the resolution property.

As above, let $Y_r$ be the smooth non-separated scheme $A^r\cup_{A^r-0}
A^r$ over a field $k$ of characteristic not 2.
Let $X_r$ be the quotient stack of $Y_r$ by the action of $\Z/2$
which is the identity outside the origin and which switches the two origins.
We will show that $X_r$ has the desired properties for $r\geq 2$.

For $r=1$, $X_1$ does have the resolution property. Indeed, $Y_1$
is a smooth scheme with affine diagonal and hence has an ample family
of line bundles by Brenner and Schr\"oer, as mentioned in section
\ref{history}. Therefore the quotient stack $X_1=Y_1/\Z/2$ has the
resolution property by Theorem \ref{thom}. More explicitly,
the stack $X_1$ has the resolution property, by Theorem \ref{thom}, because
it is the quotient of the quasi-affine scheme $A^2-0$ by the orthogonal
group $O(2)$. (The orthogonal group $O(2)$ of the quadratic form $x_1x_2$
is a split extension of $\Z/2$ by $G_m$.)

For $r\geq 2$, one checks (by arguments as in Example 1) that
pulling back via the flat morphism $X_r\arrow A^r/\Z/2$ induces
an equivalence of categories of vector bundles. Here $A^r/\Z/2$ denotes
the quotient stack of $A^r$ by the trivial action of $\Z/2$, and so a vector
bundle on $A^r/\Z/2$
is simply a direct sum of bundles $E_0\oplus E_1$ on $A^r$, where $\Z/2$
acts trivially on $E_0$ and by $-1$ on $E_1$. It follows that
\begin{align*}
K_0^{\naive}X &=K_0^{\naive}(A^r/\Z/2)\\
&= \Rep(\Z/2)\\
&= \Z^2.
\end{align*}

The true $K$-group $K_0X_r$ maps isomorphically to $G_0X_r$ since
$X_r$ is smooth over $k$. By the localization sequence as in Example 1,
$K_0X_r$ is the quotient of $K_0^{\naive}X_r=\Rep(\Z/2)$ by the relation
that
$$W\cdot \lambda_{-1}V=(\text{dim }W)\cdot \lambda_{-1}V$$
for all $W\in\Rep(\Z/2)$, where $V$ is the representation of $\Z/2$
on $A^r$. In this example, $V$ is the trivial representation, and so
$\lambda_{-1}V=0$. Therefore
the map from
$K_0^{\naive}X_r$ to $K_0X_r$ is an isomorphism, as promised.

Finally, we know that $X_r$ does not have the resolution property by
Proposition
\ref{affine}, since $X_r$ does not have affine diagonal, using that $r\geq
2$. One can define an explicit coherent sheaf $K$ on $X_r$ which is not
a quotient of a vector bundle, by the same formula as in Example 1.

\section{How to prove the resolution property in an example: the nodal cubic}
\label{nodal}

Suppose that one wishes to prove the resolution property for a stack $X$.
The proof of Theorem \ref{main} gives an idea of how to proceed. 
In many cases, the proof indicates how to construct a coherent sheaf $C$ on
$X$ such that $X$ has the resolution property if and only if the single
sheaf
$C$ is a quotient of a vector bundle. One statement of this type
is formulated in Lemma \ref{cover}, below.
I hope that this will be a useful way
to prove the resolution property in cases of interest. 

In this section, we carry the procedure out in the following
example.

\begin{proposition}
\label{propnodal}
Let $X$ be the nodal cubic over a field $k$,
that is, $\P^1$ with
the points 0 and $\infty$ identified. Let $T:=G_m$ act
on $X$ in the natural way. Then the quotient stack $X/T$
has the resolution property.
\end{proposition}

Here $X$ is a projective variety and so
the resolution property is well known for $X$, but the action of $T$
is  ``bad'' in several ways. In particular, any $T$-equivariant line
bundle on $X$ has degree 0 and so there is no $T$-equivariant embedding
of $X$ into projective space; there is not even an ample family of
$T$-equivariant
line bundles on $X$.
 Thus Theorem \ref{thom}, due to Thomason,
does not immediately apply
 to show that $X/T$ has the resolution property. The phenomenon that
$X$ has an ample family of line bundles but no ample family of
$T$-equivariant
line bundles can only occur for non-normal schemes such as this one,
which is why it seemed worth finding out whether $X/T$ has the resolution
property. (Since $X$ is not
normal, Theorem \ref{main} as stated does not apply to $X$, but the methods
still work.) The fact that we will prove the resolution property in this
``bad'' case
is encouraging for Question 3 in the introduction,
proposing that the resolution property is always preserved upon taking
the quotient by a linear algebraic group.

{\bf Proof of Proposition \ref{propnodal}. }It seems convenient
to begin by considering an etale double covering $Y$ of
$X$,
the union of two copies $A$ and $B$ of $\P^1$, with 0 in $A$ identified
with $\infty$ in $B$ and $\infty$ in $A$ identified with $0$ in $B$. I will
write
$p$ for the point $0$ in $A\subset Y$ and $q$ for the point $\infty$ in
$A\subset Y$.
The $T$-action
on $X$ lifts to a $T$-action on $Y$ which commutes with the $\Z/2$-action
(switching the two copies of $\P^1$ in the natural way), and so we have
an isomorphism of quotient stacks $X/T\cong Y/(\Z/2\times T)$. The scheme
$Y$ resembles $X$ in that it has no ample family of $T$-equivariant line
bundles. What suggests that $Y$ should be easier to study
than $X$ is that
unlike $X$, $Y$ is at least a union of $T$-invariant affine open subsets,
$Y-p$ and $Y-q$.

To show that the stack $Y/T$ has the resolution property, we will
use the following lemma, which isolates part of the proof of Theorem
\ref{main}.

\begin{lemma}
\label{cover}
Let $Y$ be a noetherian scheme with an action of a flat affine
group scheme $T$ of finite type over $\Z$ or over a field.
Let $S_1,\ldots,S_r$ be closed $T$-invariant
subschemes whose complements form an
affine open covering of $Y$.  Let $C$ be the direct sum of the ideal
sheaves $I_{S_1},\ldots,I_{S_r}$. Then the stack $Y/T$ has the
resolution property if and only if the $T$-equivariant coherent
sheaf $C$ on $Y$ is a quotient of some $T$-equivariant vector bundle
on $Y$.
\end{lemma}

{\bf Proof. }To say that the stack $Y/T$ has the resolution property
means that every $T$-equivariant coherent sheaf on $Y$ is a quotient
of some $T$-equivariant vector bundle on $Y$.
So suppose that $C$ is a quotient
of some $T$-equivariant vector bundle $E$ on $Y$.  Let $W$ be the
$GL(n)$-bundle
over $Y$ corresponding to $E$. Then the pullback of $C$ to $W$ is spanned
by its global sections. By the choice of $C$, plus affineness of the
morphism
$W\arrow Y$, it follows that the scheme $W$ is quasi-affine (this
argument is given in more detail in the proof of Theorem \ref{main}).
Therefore the stack
$W/(T\times GL(n))\cong  Y/T$ has the resolution property by
Theorem \ref{thom}. \qed

We return to the union $Y$ of two copies of $\P^1$ with the action
of $T=G_m$. 
Since $Y$ is the union of the $T$-invariant affine open subsets $Y-p$
and $Y-q$, Lemma \ref{cover} shows that $Y$ has the resolution property
if the $T$-equivariant coherent sheaf $I_{S_1}\oplus I_{S_2}$ on $Y$ is 
a quotient of some $T$-equivariant vector bundle, for some $T$-invariant
subschemes $S_1$ and $S_2$ with support equal to $p$ and $q$, respectively.
By the $\Z/2$-symmetry of $Y$, it suffices to show that
$I_S$ is a quotient of a $T$-equivariant vector bundle on $Y$,
for some $T$-invariant subscheme $S$ with support equal to $p$. Replacing
the bundle by its dual, it is equivalent to find a $T$-equivariant vector
bundle
$\E $ on $Y$ with a $T$-invariant section $s:O_Y\arrow \E $ which vanishes
(to any order) at $p$ and nowhere else on $Y$.

To define a $T$-equivariant vector bundle $\E $ on $Y$, we need
to define $T$-equivariant vector bundles $E$ on $A$ and $F$ on $B$,
together with $T$-equivariant isomorphisms $E|_0\cong F|_{\infty}$
and $E|_{\infty}\cong F|_0$. Take $E=O(1)\oplus O(-1)$ and
$F=O(1)\oplus O(-1)$ as vector bundles on $\P^1$. Define the action of $T$
on $E$ to be trivial on $E|_{\infty}$, with weight 1 on $O(1)|_0\subset
E|_0$, and weight $-1$ on $O(-1)|_{\infty}\subset E|_0$. Define the action
of $T$ on $F$ to be trivial on $F|_0$, with weight $-1$ on
$O(1)|_{\infty}\subset F|_{\infty}$, and weight $1$ on $O(-1)|_{\infty}
\subset F|_{\infty}$. Clearly there are isomorphisms of
$T$-representations $E|_0\cong F|_{\infty}$ and
$E|_{\infty}\cong F|_0$, which we can use to define a $T$-equivariant
vector bundle $\E =(E,F)$ on $Y$. For our purpose, we need to choose
the isomorphism
between the trivial 2-dimensional representations $E|_{\infty}$
and $F|_0$ of $T$ so as to send the line $O(1)|_{\infty}
\subset E|_{\infty}$ to the line $O(1)|_0\subset F|_0$. With this choice,
the $T$-vector bundle $\E $ on $Y$
is not a direct sum of two $T$-line bundles.

The vector bundle $\E $ has a $T$-invariant section over $Y$ (contained
in the subbundle $O(1)\subset E$ over $A$ and in the subbundle
$O(1)\subset F$ over $B$) which vanishes at $p$ (the image of
$0\in A\cong \P^1$ and of $\infty\in B\cong \P^1$) but nowhere else in $Y$.
Thus, as we have explained, the dual $T$-vector bundle $\E ^*$ on $Y$
maps onto the $T$-coherent sheaf $I_S$ for some subscheme $S$
with support equal to $p$. By Lemma \ref{cover},
the stack $Y/T$ has the resolution property.

From here it is easy to deduce that the quotient stack of
the nodal cubic $X$ by $T$ also has the resolution property,
as we wanted. Namely,
writing $\sigma$ for the generator of the $\Z/2$-action on $Y$,
$\E ^*\oplus \sigma^*(\E ^*)$ is a $\Z/2\times T$-equivariant vector bundle
on $Y$ which maps onto the $\Z/2\times T$-equivariant sheaf $I_S\oplus
I_{\sigma(S)}$. So, by Lemma \ref{cover} again, the stack $Y/(\Z/2
\times T)=X/T$ has the resolution property. Proposition
\ref{propnodal} is proved. \qed

It is interesting to add that the proof of Theorem \ref{main} works
completely
for $X/T$ even though the nodal cubic
$X$ is not normal, because it is a scheme, not just an algebraic space.
That proof shows that since $X/T$ has the resolution property, it is the
quotient of a quasi-affine scheme by $GL(n)$ for some $n$. Thus this
apparently
bad example of a stack
turns out to have some very good properties. As mentioned at the beginning
of this section, this is encouraging
for Question 3 in the introduction.


\small \sc DPMMS, Wilberforce Road,
Cambridge CB3 0WB, England.

b.totaro@dpmms.cam.ac.uk
\end{document}